# Two Applications of Desargues' Theorem

Prof. Florentin Smarandache, University of New Mexico, U.S.A.
Prof. Ion Pătrașcu, The National College "Frații Buzești", Craiova, Romania

In this article we will use the Desargues' theorem and its reciprocal to solve two problems.

For beginning we will enunciate and prove Desargues' theorem:

**Theorem 1** (G.Desargues, 1636, the famous "perspective theorem": When two triangles are in perspective, the points where the corresponding sides meet are collinear.)

Let two triangle $ABC$ and $A_1B_1C_1$ be in a plane such that $AA_1 \cap BB_1 \cap CC_1 = \{O\}$,

$$AB \cap A_1B_1 = \{N\}$$
$$BC \cap B_1C_1 = \{M\}$$
$$CA \cap C_1A_1 = \{P\}$$

then the points $N$, $M$, $P$ are collinear.

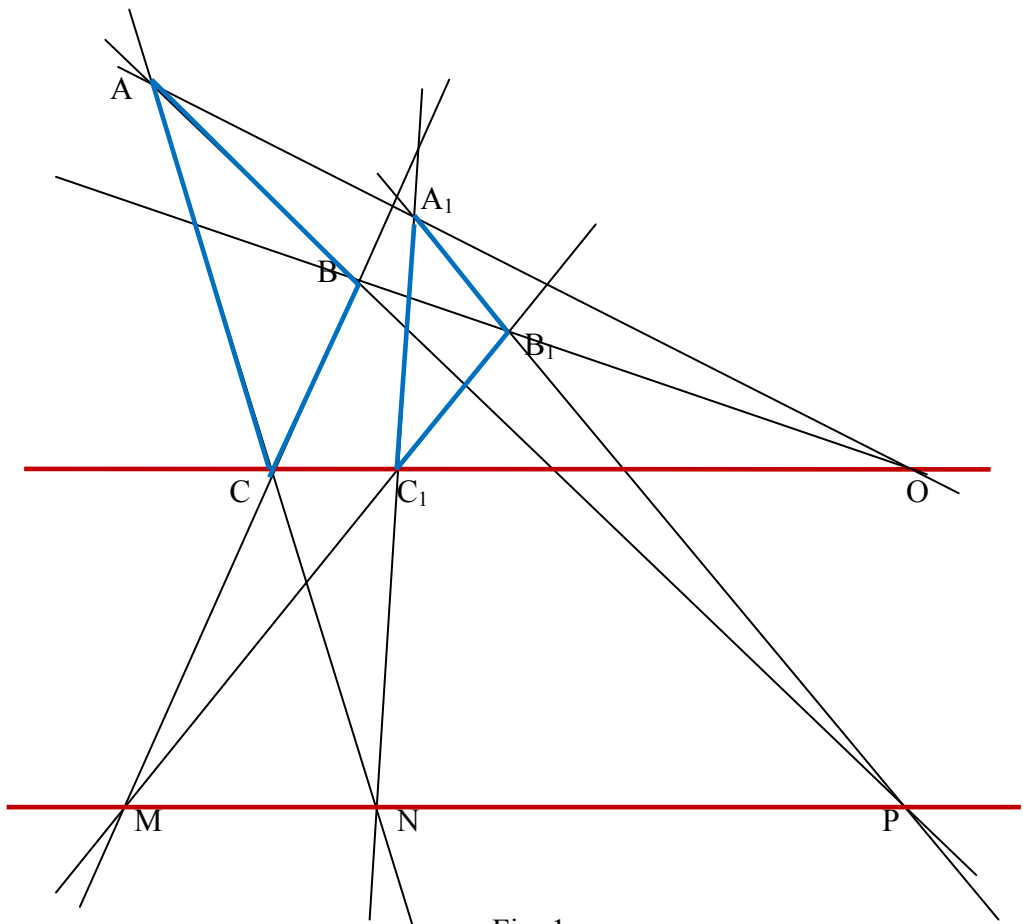

Fig. 1

**Proof**

Let $\{O\} = AA_1 \cap BB_1 \cap CC_1$, see Fig.1.. We'll apply the Menelaus' theorem in the triangles $OAC$; $OBC$; $OAB$ for the transversals $N, A_1, C_1$; $M, B_1, C_1$; $P, B_1, A_1$, and we obtain

$$\frac{NA}{NC} \cdot \frac{C_1C}{C_1O} \cdot \frac{A_1O}{A_1A} = 1 \qquad (1)$$

$$\frac{MC}{MB} \cdot \frac{B_1B}{B_1O} \cdot \frac{C_1O}{C_1C} = 1 \qquad (2)$$

$$\frac{PB}{PA} \cdot \frac{B_1O}{B_1B} \cdot \frac{A_1A}{A_1O} = 1 \qquad (3)$$

By multiplying the relations (1), (2), and (3) side by side we obtain

$$\frac{NA}{NC} \cdot \frac{MC}{MB} \cdot \frac{PB}{PA} = 1.$$

This relation, shows that $N, M, P$ are collinear (in accordance to the Menealaus' theorem in the triangle $ABC$).

**Remark 1**

The triangles $ABC$ and $A_1B_1C_1$ with the property that $AA_1, BB_1, CC_1$ are concurrent are called homological triangles. The point of concurrency point is called the homological point of the triangles. The line constructed through the intersection points of the homological sides in the homological triangles is called the triangles' axes of homology.

**Theorem 2** (The reciprocal of the Desargues' theorem)

` If two triangles $ABC$ and $A_1B_1C_1$ are such that

$$AB \cap A_1B_1 = \{N\}$$
$$BC \cap B_1C_1 = \{M\}$$
$$CA \cap C_1A_1 = \{P\}$$

And the points $N, M, P$ are collinear, then the triangles $ABC$ and $A_1B_1C_1$ are homological.

**Proof**

We'll use the reduction ad absurdum method.

Let

$$AA_1 \cap BB_1 = \{O\}$$
$$AA_1 \cap CC_1 = \{O_1\}$$
$$BB_1 \cap CC_1 = \{O_2\}$$

We suppose that $O \neq O_1 \neq O_2 \neq O_3$.

The Menelaus' theorem applied in the triangles $OAB$, $O_1AC$, $O_2BC$ for the transversals $N, A_1, B_1$; $P, A_1, C_1$; $M, B_1, C_1$, gives us the relations

$$\frac{NB}{NA} \cdot \frac{B_1O}{B_1B} \cdot \frac{AA_1}{A_1O} = 1 \qquad (4)$$

$$\frac{PA}{PC} \cdot \frac{A_1O_1}{A_1O} \cdot \frac{C_1C}{C_1O_1} = 1 \qquad (5)$$

$$\frac{MC}{MB} \cdot \frac{B_1B}{B_1O} \cdot \frac{C_1O_2}{C_1C} = 1 \qquad (6)$$

Multiplying the relations (4), (5), and (6) side by side, and taking into account that the points $N$, $M$, $P$ are collinear, therefore

$$\frac{PA}{PC} \cdot \frac{MC}{MB} \cdot \frac{NB}{NA} = 1 \qquad (7)$$

We obtain that

$$\frac{A_1O_1}{A_1O} \cdot \frac{B_1O}{B_1O_2} \cdot \frac{C_1O_2}{C_1O_2} = 1 \qquad (8)$$

The relation (8) relative to the triangle $A_1B_1C_1$ shows, in conformity with Menelaus' theorem, that the points $O, O_1, O_2$ are collinear. On the other hand the points $O, O_1$ belong to the line $AA_1$, it results that $O_2$ belongs to the line $AA_1$. Because $BB_1 \cap CC_1 = \{O_2\}$, it results that $\{O_2\} = AA_1 \cap BB_1 \cap CC_1$, and therefore $O_2 = O_1 = O$, which contradicts the initial supposition.

**Remark 2**
The Desargues' theorem is also known as the theorem of the homological triangles.

**Problem 1**
If $ABCD$ is a parallelogram, $A_1 \in (AB), B_1 \in (BC), C_1 \in (CD), D_1 \in (DA)$ such that the lines $A_1D_1, BD, B_1C_1$ are concurrent, then:

a) The lines $AC, A_1C_1$ and $B_1D_1$ are concurrent

b) The lines $A_1B_1, C_1D_1$ and $AC$ are concurrent.

**Solution**

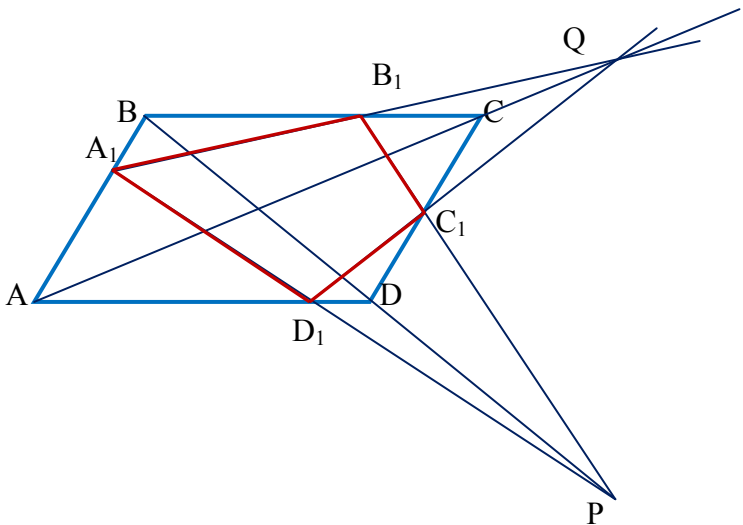

Fig. 2

Let $\{P\} = A_1D_1 \cap B_1C_1 \cap BD$ see Fig. 2. We observe that the sides $A_1D_1$ and $B_1C_1$; $CC_1$ and $AD_1$; $A_1A$ and $CB_1$ of triangles $AA_1D_1$ and $CB_1C_1$ intersect in the collinear points $P, B, D$. Applying the reciprocal theorem of Desargues it results that these triangles are homological, that is, the lines: $AC, A_1C_1$ and $B_1D_1$ are collinear.

Because $\{P\} = A_1D_1 \cap B_1C_1 \cap BD$ it results that the triangles $DC_1D_1$ and $BB_1A_1$ are homological. From the theorem of the of homological triangles we obtain that the homological lines
$DC_1$ and $BB_1$; $DD_1$ and $BA_1$; $D_1C_1$ and $A_1B_1$ intersect in three collinear points, these are $C, A, Q$, where $\{Q\} = D_1C_1 \cap A_1B_1$. Because $Q$ is situated on $AC$ it results that $A_1B_1, C_1D_1$ and $AC$ are collinear.

**Problem 2**

Let $ABCD$ a convex quadrilateral such that
$$AB \cap CB = \{E\}$$
$$BC \cap AD = \{F\}$$
$$BD \cap EF = \{P\}$$
$$AC \cap EF = \{R\}$$
$$AC \cap BD = \{O\}$$

We note with $G$, $H$, $I$, $J$, $K$, $L$, $M$, $N$, $Q$, $U$, $V$, $T$ respectively the middle points of the segments: $(AB), (BF), (AF), (AD), (AE), (DE), (CE), (BE), (BC), (CF), (DF), (DC)$. Prove that

i) The triangle $POR$ is homological with each of the triangles: $GHI$, $JKL$, $MNQ$, $UVT$.

ii) The triangles $GHI$ and $JKL$ are homological.

iii) The triangles $MNQ$ and $UVT$ are homological.

iv) The homology centers of the triangles $GHI$, $JKL$, $POR$ are collinear.

v) The homology centers of the triangles $MNQ$, $UVT$, $POR$ are collinear.

**Solution**

i) when proving this problem we must observe that the $ABCDEF$ is a complete quadrilateral and if $O_1, O_2, O_3$ are the middle of the diagonals $(AC), (BD)$ respective $EF$, these point are collinear. The line on which the points $O_1, O_2, O_3$ are located is called the Newton-Gauss line [* for complete quadrilateral see [1]].

The considering the triangles $POR$ and $GHI$ we observe that $GI \cap OR = \{O_1\}$ because $GI$ is the middle line in the triangle $ABF$ and then it contains the also the middle of the segment $(AC)$, which is $O_1$. Then $HI \cap PR = \{O_3\}$ because $HI$ is middle line in the triangle $AFB$ and $O_3$ is evidently on the line $PR$ also. $GH \cap PO = \{O_2\}$ because $GH$ is middle line in the triangle $BAF$ and then it contains also $O_2$ the middle of the segment $(BD)$.

The triangles $GIH$ and $ORP$ have as intersections of the homological lines the collinear points $O_1, O_2, O_3$, according to the reciprocal theorem of Desargues these are homological.

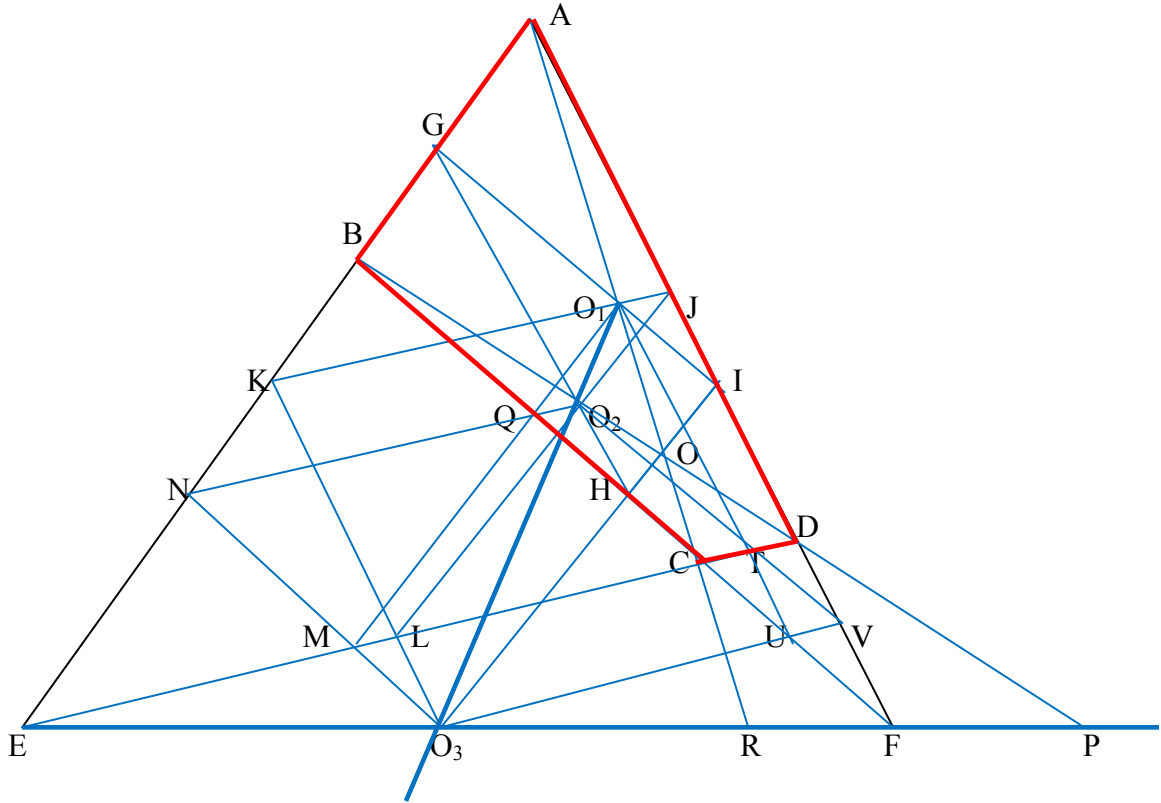

Fig. 3

Similarly, we can show that the triangle $ORP$ is homological with the triangles $JKL$, $MNQ$, and $UVT$ (the homology axes will be $O_1, O_2, O_3$).

ii)     We observe that
$$GI \cap JK = \{O_1\}$$
$$GH \cap JL = \{O_2\}$$
$$HI \cap KL = \{O_3\}$$
then $O_1, O_2, O_3$ are collinear and we obtain that the triangles $GIH$ and $JKL$ are homological

iii)     Analog with ii)

iv)     Apply the Desargues' theorem. If three triangles are homological two by two, and have the same homological axes then their homological centers are collinear.

v)     Similarly with iv).

**Remark 3**

The precedent problem could be formulates as follows:

The four medial triangles of the four triangles determined by the three sides of a given complete quadrilateral are, each of them, homological with the diagonal triangle of the complete quadrilateral and have as a common homological axes the Newton-Gauss line of the complete quadrilateral.

We mention that:
- The *medial triangle* of a given triangle is the triangle determined by the middle points of the sides of the given triangle (it is also known as the complementary triangle).
- The *diagonal triangle* of a complete quadrilateral is the triangle determined by the diagonals of the complete quadrilateral.

We could add the following comment:

Considering the four medial triangles of the four triangles determined by the three sides of a complete quadrilateral, and the diagonal triangle of the complete quadrilateral, we could select only two triplets of triangles homological two by two. Each triplet contains the diagonal triangle of the quadrilateral, and the triplets have the same homological axes, namely the Newton-Gauss line of the complete quadrilateral.

**Open problems**
1. What is the relation between the lines that contain the homology centers of the homological triangles' triplets defined above?
2. Desargues theorem was generalized in [2] in the following way: Let's consider the points $A_1,...,A_n$ situated on the same plane, and $B_1,...,B_n$ situated on another plane, such that the lines $A_iB_i$ are concurrent. Then if the lines $A_iA_j$ and $B_iB_j$ are concurrent, then their intersecting points are collinear.
   Is it possible to generalize Desargues Theorem for two polygons both in the same plane?
3. What about Desargues Theorem for polyhedrons?